\documentclass[12pt,twoside]{paper}
\usepackage{amsmath}
\usepackage{hyperref}
\usepackage{amsfonts}
\usepackage{amssymb}
\usepackage{amsthm}

\usepackage{graphicx}

\def\bee{\begin{equation}}
\def\eee{\end{equation}}
\def\bdm{\begin{displaymath}}
\def\edm{\end{displaymath}}

\def\Mas{Ma$\acute{\rm s}$lanka  }
\def\BD{B{\'a}ez-Duarte  }
\begin{document}

\thispagestyle{empty}
\centerline{}
\bigskip
\bigskip
\bigskip
\bigskip
\bigskip
\bigskip
\centerline{\Large\bf Failed attempt to disproof  the Riemann Hypothesis}
\bigskip

\begin{center}
{\large \sl Marek Wolf}\\*[5mm]

% Institute of Theoretical Physics, University of Wroc{\l}aw\\
%Pl.Maxa Borna 9, PL-50-204 Wroc{\l}aw, Poland,
\href{mailto:~mwolf@ift.uni.wroc.pl}{e-mail:~mwolf@ift.uni.wroc.pl}\\
\bigskip

\end{center}

\bigskip
\begin{center}
{\bf Abstract}\\
\end{center}

\begin{minipage}{12.8cm}
In this paper we are going to describe the results of the computer experiment,
which in principle can rule out the Riemann Hypothesis. We use the sequence
$c_k$ appearing in the \BD   criterion for the RH.
Namely we calculate $c_{100000}$ with thousand  digits of accuracy using two
different formulas for $c_k$ with the aim to disproof the Riemann
Hypothesis in the case these two numbers will differ. We found the discrepancy
only on the 996 decimal place (accuracy of $10^{-996}$). The reported here experiment can
be of interest for developers of Mathematica and PARI/GP.
\end{minipage}

\bigskip\bigskip

\bibliographystyle{abbrv}

\section{Introduction}

With the advent of the computers era the  computing machines have been used to prove
mathematical theorems. The most
spectacular examples of  such a use of computers were proofs of the
four color theorem  \cite{Appel_Haken1}, \cite{Appel_Haken2}
and of the Kepler conjecture about sphere
packing in three-dimensional Euclidean space \cite{Kepler}.
It seems to be not possible to use computers for the proof of the Riemann Hypothesis (RH),
but refutation of it by numerical calculations seems to be plausible.

The Riemann Hypothesis says that  the series
\bee
\zeta(s)=\sum_{n=1}^\infty \frac{1}{n^s}, ~~~~~~(s=\sigma+it, ~~~~\Re\{s\}>1)
\label{zeta}
\eee
analytically continued to the complex plane in addition to trivial zeros $\zeta(-2n)=0$
has nontrivial zeros  $\zeta(\rho_l)=0$
in the  critical strip $0<\Re\{s\}<1$ only on the critical line:
$\Re\{\rho_l\}=\frac{1}{2}$  i.e. $\rho_l=\frac{1}{2}+i\gamma_l$, see e.g. the modern
guide to the RH \cite{Borwein_RH}. In the same book there is a review
of failed attempts to prove RH in Chapter 8.  Presently the requirement that the
nontrivial zeros are  simple  $\zeta'(\rho_l)\neq 0$  is often added.

The first use of computers in connection with RH was  checking by Allan Turing,
whether the nontrivial zeros of $\zeta(s)$ have indeed real part $\frac{1}{2}$
\cite{Turing1953}.  Turing was convinced that RH is false. Let us quote the
sentence  from the first page of his paper: ``The calculation were done in an
optimistic hope that a zero would be found
off the critical line'', but up to $t=1540$ Turing has found that all zeros are on
the critical line. The present record belongs to Xavier Gourdon \cite{Gourdon}
who has checked that all $10^{13}$ first zeros of the Riemann $\zeta(s)$ lie on
the critical line. Andrew Odlyzko checked that RH is true in different intervals
around $10^{20}$ \cite{Odlyzko_zera1}, $10^{21}$ \cite{Odlyzko_zera2}, $10^{22}$
\cite{Odlyzko2001}, but his aim was not verifying the RH but rather providing evidence
for conjectures that relate nontrivial zeros  of $\zeta(s)$ to eigenvalues of
random matrices. In fact Odlyzko has expressed the view, that off critical line zeros
could be encountered at least at $t$ of the order $10^{10^{10000}}$, see
\cite[p.358]{Derbyshire}. Asked by Derbyshire ``What do you think about this
darn Hypothesis? Is it true, or not?'' Odlyzko replied:  ``Either it's true, or else
it isn't''.  Also other famous mathematicians John E. Littlewood and
Paul Turan have not believed RH is true.

There were several attempts  to use computers to disprove some conjectures
related to RH
in the past.  Sometimes it was  sufficient to find a counterexample to  the given
hypothesis, sometimes the disprove was not direct. For example, the
Haselgrove \cite{Haselgrove} disproved the P\'olya's Conjecture stating that the
function
\bee
L(x):=\sum_{n\leq x}\lambda(n)
\eee
satisfies $L(x)\leq 0$ for $x\geq 2$, where $\lambda(n)$ is Liouville's function
defined by
\bdm
\lambda(n) =(-1)^{r(n)}
\edm
where $r(n)$ is the number of, not necessarily distinct, prime factors in
$n=p_1^{r_1}\cdots p_{\alpha(n)}^{r_n}$, with multiple factors
counted multiply: $r=r_1+\ldots+r_n$.  %,  $\lambda(1)=0, \lambda(12)=-1.$
From the truth of the Polya Conjectures
the RH follows, but not the other way around.  The Haselgrove proof was indirect,
and in 1960 Lehman \cite{Lehman1960} has found on the computer explicit
counterexample: $ L(906180359)=1$.

The next example is provided by the Mertens conjecture. Let $M(x)$ denote
the Mertens function defined by
\bee
M(x)=\sum_{n<x} \mu(n), 	
\eee
where $\mu(n)$ is the M{\"o}bius function
\bdm
\mu(n) \,=\,
\left\{
\begin{array}{ll}
1 & \mbox {for  $ n =1 $} \\
0 & \mbox {when $p^2|n$}\\
(-1)^r & \mbox{\rm when}~ n=p_1 p_2 \ldots p_r
\end{array}
\right.
\edm
From
\bee
|M(x)|<x^\frac{1}{2}
\eee
again the RH would follow. However in 1985 A. Odlyzko and H. te Riele \cite{odlyzko1985}
disproved  the Mertens conjecture, again not directly, but  later it was shown
by J. Pintz \cite{Pintz1987}  that
the first counterexample appears below $exp(3.21 \times 10^{64})$.
The upper bound has since been
lowered to $exp(1.59\times 10^{40})$ \cite{KotnikRiele2006}.

Especially  interesting is the value of the de Bruijn-Newman constant $\Lambda$,
see e.g. \S2.32 (pp. 203-205) in \cite{Finch}. Unconditionally  $\Lambda\leq 1/2$
and the Riemann Hypothesis is equivalent
to the inequality $\Lambda\leq 0$. The fascinating run for the best lower bound on
$\Lambda$  ended with the value $\Lambda>-2.7\times 10^{-9}$ obtained by  Odlyzko
\cite{Odlyzko2000}. Such a narrow gap for values of  $\Lambda$ being compatible
with RH allowed Odlyzko to make the remark: ``the Riemann Hypothesis, if true,
is just barely true''.

In 1997 Xian-Jin  Li proved \cite{Li-1997}, that  Riemann Hypothesis is true  {\it iff}
the sequence:
\begin{displaymath}
\lambda_n=\frac{1}{(n-1)!}\frac{d^n}{ds^n}(s^{n-1}\log \xi(s)) |_{s=1} %\vert_{s=1}
\end{displaymath}
where
\begin{displaymath}
\xi(s)=\frac{1}{2}s(s-1)\Gamma\left(\frac{s}{2}\right)\zeta(s)
\end{displaymath}
fulfills:
\bee
\lambda_n\ge 0~~~~ {\rm for}~~  n=1,2,\dots
\label{Li}
\eee
The explicit expression has the form:
\bee
~~\lambda_n =\sum_\rho (1-(1-1/\rho)^n).
\eee
K. Ma{\'s}lanka \cite{Maslanka-2004}, \cite{Maslanka-2004-b}
has performed extensive computer calculations of these constants confirming (\ref{Li}).

Let us mention also  the elementary
Lagarias criterion \cite{Lagarias}: to disprove the RH it suffices to find one $n$
that has so  many divisors, that:
\bee
 \sum_{d|n} d > H_n + \exp(H_n)\log(H_n).
\eee
The  Lagarias criterion is not well suited for computer verification (it is not an easy
task to calculate  $H_n$ for $n \sim 10^{100000}$ with sufficient accuracy)
and in \cite{Briggs2006}
Keith Briggs has undertaken instead the verification of the Robin \cite{Robin}
criterion for RH:
\bee
{\rm RH~~} \Leftrightarrow  \sum_{d|n} d < e^\gamma n \log\log(n)~~~~~~~~~{\rm for~~ }n>5040
\eee
For appropriately chosen $n$ Briggs obtained for the difference between
r.h.s. and l.h.s. of the above inequality value as small as
$e^{-13}\approx 2.2\times 10^{-6}$, hence again RH is very close to be violated.

In this paper we are going to propose a method  which in principle can provide
a refutation of the RH. The idea is to calculate some number with very high accuracy
(one  thousand  digits) in two ways: one without any knowledge on the zeros of $\zeta(s)$
and second using the explicit formula involving   all $\rho_l$.
Despite  some estimation presented in Sect.3 indicating that
the discrepancy could be found merely at much higher than a thousand decimal place we
performed the calculations in an  optimistic hope that we will find the discrepancy
between these two numbers, paraphrasing the sentence of Turing.  There is
a lot of number theoretic functions defined often  in an elementary way
being expressed also by the ``explicit'' formulas in terms of zeros
of the $\zeta(s)$ function. Let us mention here  the Chebyshev function
\[
\psi(x) = \sum_{n<x} \Lambda(n) ~~(=\log({\rm lcm}(2, 3, \cdots \lfloor x \rfloor))),
\]
where von  Mangoldt function $\Lambda(n)$ is defined as
\bdm
\Lambda(n) \,=\,
\left\{
\begin{array}{ll}
\log p & \mbox {for $ n =p^m $} \\
0 & \mbox {in other cases }
\end{array}
\right.
\edm
The explicit formula reads, see eg. \cite{Titchmarsh}:
\bee
\psi(x) =  x -\frac{\zeta'(0)}{\zeta(0)}  - \!\!\!\!\!\!\sum_{all~zeros~\rho_l} \!\!\frac{x^{\rho_l}}{\rho_l} =
x - \log(2\pi)   - \frac{1}{2}\log\left(1-\frac{1}{x^2}\right) ~ - \!\!\!\!\!\! \sum_{nontr. ~zeros ~\rho_l}\!\!\! \frac{x^{\rho_l}}{\rho_l}
\label{psi-explicite}
\eee

Also the Mertens's function has the explicit representation (last term is comprising
contribution from all trivial zeros) \cite{Titchmarsh}:
\bee
\sum_{n<x} \mu(n) = \sum_{nontr. ~zeros~ \rho_l} \frac{x^{\rho_l}}{\rho_l\zeta'(\rho_l)} -
2 -\sum_{n=1}^\infty (-1)^n \left(\frac{2\pi}{x}\right)^{2n}\frac{1}{(2n)!n\zeta(2n+1)}
\eee
The problem with these series is they are extremely slow to converge because the partial
sums oscillate with amplitudes diminishing at very slow rates. For example $\psi(1000001)=
999586.597\ldots$, while from (\ref{psi-explicite}) summing over 5,549,728 zeros gives
$999587.15\ldots$, thus relative error is 0.000055.

In the  computer  experiment reported  here we were able to get discrepancy less than
$10^{-996}$ between  quantity calculated from generic formula and from explicit one
summed over only 2600 nontrivial zeros computed with 1000 significant digits.

\section{The \BD criterion for the Riemann Hypothesis}

We begin by recalling the
following representation of the $\zeta(s)$ function valid on the whole
complex plane without $s=1$ found by Krzysztof \Mas \cite{Maslanka}:
\begin{equation}
\zeta (s)=\frac{1}{s-1}\sum_{k=0}^{\infty }\frac{\Gamma \left( k+1-%
\frac{s}{2}\right) }{\Gamma \left( 1-\frac{s}{2}\right) }\frac{A_{k}}{k!}%
\label{zeta-Maslanka}
\end{equation}%
where%
\begin{equation}
A_{k}:=\sum_{j=0}^{k}(-1)^{j}\binom{k}{j}(2j+1)\zeta (2j+2)\equiv
\sum_{j=0}^{k}\binom{k}{j}\frac{\pi ^{2j+2}}{\left( 2\right)
_{j}\left( \frac{1}{2}\right) _{j}}B_{2j+2}.  
\label{ak}
\end{equation}

The expansion (\ref{zeta-Maslanka}) provides
an example of the analytical continuation of (\ref{zeta}) to the  whole complex plane
except $s=1$.  Since $A_{k}$ tend to zero fast  as
$k\rightarrow \infty $ the expansion (\ref{zeta-Maslanka}) converges uniformly on the
whole complex plane \cite{baezduarte-2003}. Based on the representation
(\ref{zeta-Maslanka}) Luis \BD in \cite{Baez-Duarte} has proved that
RH is equivalent to the statement that
\begin{equation}
c_{k} = \mathcal{O}\left(  k^{-3/4+\varepsilon }\right),~~~~\forall \varepsilon >0
\label{kryterium}
\end{equation}
where
\bee
c_{k} = \sum_{j=0}^{k}(-1)^{j}\binom{k}{j}\frac{1}{\zeta (2j+2)}.
\label{ck-generic}
\eee
If additionally
\bee
c_{k} = \mathcal{O}\left(  k^{-3/4}\right)
\label{kryterium_2}
\eee
then all zeros of $\zeta(s)$ are simple. \BD has shown unconditionally
(regardless of validity of the RH) slower decrease $c_{k} = \mathcal{O}(k^{-1/2})$.
The plot of $c_k$  for $k=1, 2,  \ldots, 100000$ is presented in the Fig.1.

\begin{figure}[t!]
\begin{center}
\includegraphics[width=14cm,angle=0, scale=1]{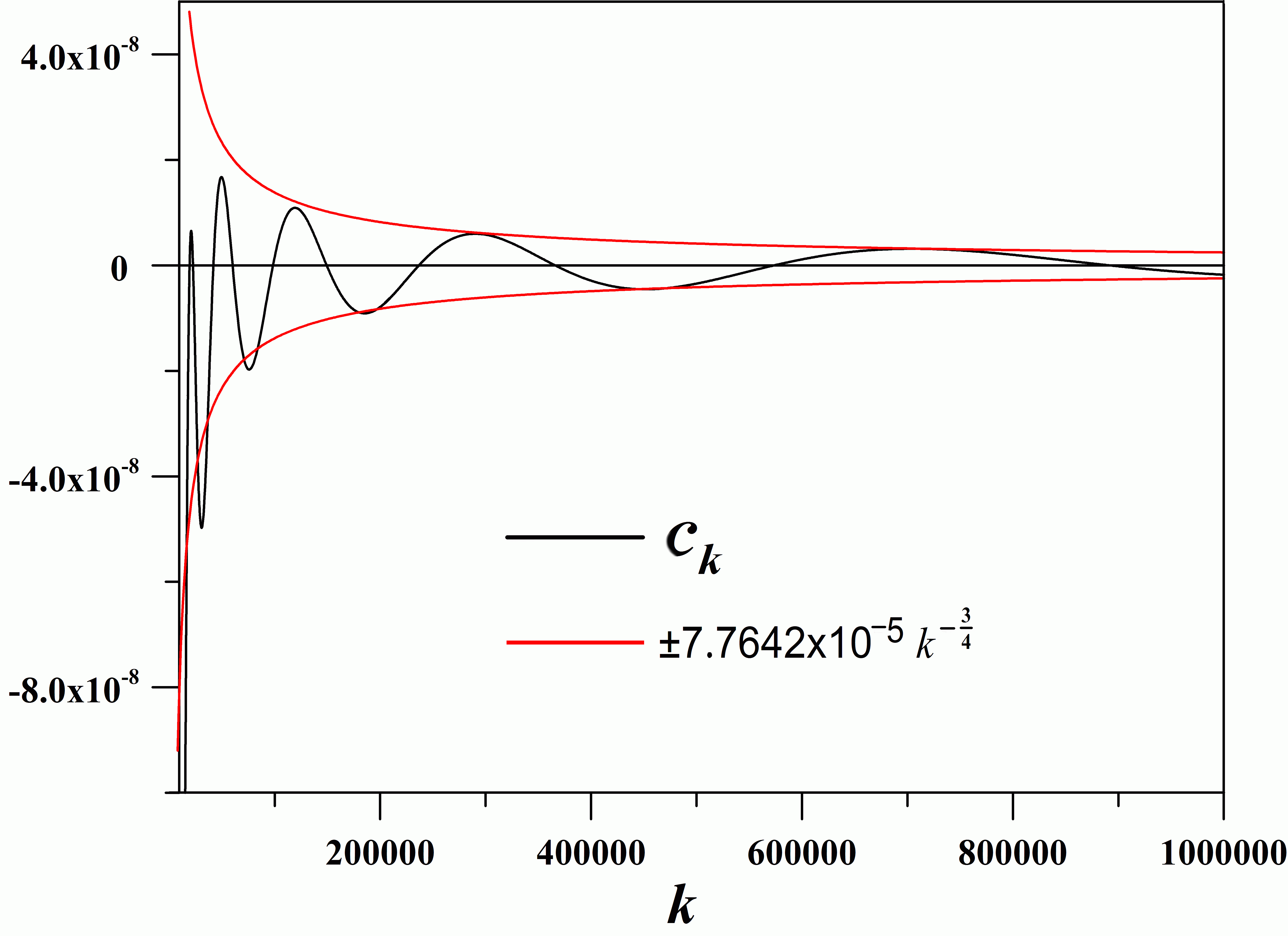} \\
Fig.1  The plot of the \BD sequence $c_k$ for $k\in(1,10^6)$. The equation
for envelope was obtained from the explicit equation (\ref{ck-oscillation}): for large
$k$ the oscillating part $\tilde{c_k}$ is dominant and for $k>100000~~c_k$ fits well
between the red lines,  for details  see  \cite{WolfCMIST}.%\\
\end{center}
\end{figure}

The explicit formula expressing the \BD sequence  $c_k$ directly in terms of the zeros
of $\zeta(s)$  can be
written as a sum of two parts: quickly decreasing with $k$ trend $\bar{c}_k$
arising from trivial zeros of $\zeta(s)$
and oscillations $\tilde{c}_k$ involving complex nontrivial zeros:
\bee
c_k=\bar{c}_k +\tilde{c}_k
\label{ck-explicite}
\eee

Because the derivatives $\zeta'(-2n)$ at trivial zeros are known analytically:
\bee
\zeta'(-2n)=\frac{(-1)^n\zeta(2n+1)(2n)!}{2^{2n+1}\pi^{2n}}.
\eee
\Mas in \cite{Maslanka3} was able to  give the closed expression for trend stemming from
zeros $\rho_n=-2n$:
\bee
\bar{c}_k=-\frac{1}{(2\pi)^2}\sum_2^\infty \frac{\Gamma(k+1)\Gamma(m)}{\Gamma(k+m+1)
\Gamma(2m-1)}\frac{(-1)^m (2\pi)^{2m}}{\zeta(2m-1)}
\label{ck-trend}
\eee
\BD  is skipping the trend $\bar{c_k}$ remarking
only that it is of the order $o(1)$ (Remark 1.6 in \cite{Baez-Duarte}). It is an easy
calculation to show (see  \cite{WolfCMIST}), that for large $k$
\bee
\bar{c}_k =
-\frac{(2\pi)^2}{2\zeta(3)}\frac{1}{k^2}+\mathcal{O}(k^{-3}),
\eee
thus  the dependence $c_k=\mathcal{O}(k^{-3/4})$  in  (\ref{kryterium}) is linked to
the oscillating part $\tilde{c}_k$.

For $\tilde{c}_k$ \BD gives the formula \cite{Baez-Duarte}:
\bee
\tilde{c}_{k-1}=\frac{1}{2k}\sum_{\rho_l} {\rm Res}\left(\frac{1}{\zeta'(s)P_k(s/2)};
s=\rho_l\right),
\label{c_k_Luis0}
\eee
where
\bee
P_k(s):=\prod_{r=1}^k \left(1-\frac{s}{r}\right)
\label{Pochhammer}
\eee
is the Pochhammer symbol.
Assuming zeros of $\zeta(s)$ are simple we can write:
\bee
\tilde{c}_{k-1}=\frac{1}{2k}\sum_{\rho_l} \frac{1}{\zeta'(\rho_l)P_k(\rho_l/2)}.
\label{c_k_Luis}
\eee
An appropriate order of summation over nontrivial zeros is assumed in (\ref{c_k_Luis0})
and (\ref{c_k_Luis}), see \cite[Theorem 1.5] {Baez-Duarte}.
Because
\bee
P_k(s)= \frac{(-1)^k\Gamma(s)}{\Gamma(k+1)\Gamma(s-k)}
\label{P-gammy}
\eee
collecting  in pairs $\rho_l$ and $\overline{\rho}$ we can convert (\ref{c_k_Luis})
to the form:
\bee
\tilde{c}_k=(-1)^{k+1}\Re\left(\sum_{\rho_l, ~\Im(\rho_l)>0}\frac{\Gamma(k+1)\Gamma(\frac{\rho_l}{2}-k-1)}
{\Gamma(\frac{\rho_l}{2})\zeta'(\rho_l)}\right).
\eee
We have found that the numerical calculation of $P_k(\rho_l/2)$
in PARI/GP directly from above product (\ref{Pochhammer}) is much slower than
the use of the $\Gamma(z)$ functions (\ref{P-gammy}).

\BD  proves in \cite[Lemma 2.2]{Baez-Duarte}
that
\bee
\lim_{k\rightarrow \infty} P_k(s) k^s = \frac{1}{\Gamma(1-s)}
\eee
thus for large $k$  we can replace  $k+1$ by $k$ and  transform (\ref{c_k_Luis})
in the following way:
\bee
\tilde{c}_{k}=\frac{1}{2k}\sum_{\rho_l} \frac{k^{\rho_l/2}}
{\zeta'(\rho_l)P_k(\rho_l/2)k^{\rho_l/2}}=\Re\left(\sum_{\Im(\rho_l)>0}
\frac{k^{\rho_l/2-1}\Gamma(1-\frac{\rho_l}{2})} {\zeta'(\rho_l)}\right).
\label{ck-od-k}
\eee
Now we assume the RH: $\rho_l=\frac{1}{2}+i\gamma_l$. Then we get for $\tilde{c}_{k}$
the overall factor $k^{-3/4}$ --- the dependence following from RH, see
(\ref{kryterium}) --- multiplied by oscillating terms:
\bee
\tilde{c}_{k}=\frac{1}{k^{3/4}}\Re\left(\sum_{\gamma_l>0} \frac{k^{i\gamma_l/2}\Gamma(\frac{3}{4}
-\frac{i}{2}\gamma_l)} {\zeta'(\frac{1}{2}+i\gamma_l)}\right).
\label{ck-oscillation}
\eee
Using the formula (6.1.45) from \cite{abramowitz+stegun}:
\bee
\lim_{|y|\rightarrow \infty} \frac{1}{\sqrt{2\pi}}|\Gamma(x+iy)|e^{\frac{1}{2}\pi|y|}
|y|^{\frac{1}{2}-x}=1
\eee
assuming RH we obtain for large $\gamma_l>0$:  % in upper half plane:
\bee
\left|\Gamma\left(\frac{3}{4} \mp \frac{i}{2}\gamma_l\right)\right|\approx\sqrt{2\pi}
e^{-\pi\gamma_l/4}\left(\frac{\gamma_l}{2}\right)^{\frac{1}{4}},
\label{Gamma-exp}
\eee
hence  we get exponential decrease of summands in the sum (\ref{ck-oscillation})
over nontrivial zeros  giving $\tilde{c}_{k}$ and (\ref{ck-oscillation}) is very
fast convergent.  Because of that if RH is true the sum
(\ref{ck-oscillation}) will be dominated by first zero $\gamma_1= 14.13472514\ldots$
what leads to the approximate expression (for details see \cite{WolfCMIST}):
\bee
\tilde{c}_k=\frac{A}{k^{\frac{3}{4}}}\sin\left(\phi -\frac{1}{2}\gamma_1\log(k)\right)
~~~A=7.775062\ldots\times 10^{-5},~~~\phi=2.592433\ldots\ \ .
\label{sinus}
\eee
For large $k$ the above formula (\ref{sinus})  gives very fast method of calculating
quite accurate  values of $c_k$,  orders faster than (\ref{ck-generic}).

In the following we will denote by $c_k^g$ the values  of the \BD sequence
obtained from
the generic formula (\ref{ck-generic}) and by  $c_k^e$ the values obtained
from explicit formula (\ref{ck-explicite}), i.e. in fact from (\ref{ck-trend})
and (\ref{ck-oscillation}) as no one zero off critical line is known.

\section{The scenario of violation of the Riemann Hypothesis}

The condition (\ref{kryterium_2}) means that the combination $ k^{\frac{3}{4}} c_k$
should be contained between two paralel lines $\pm C$ for all $k$, where $C$ is
the constant  hidden in big-$\mathcal{O}$ in  (\ref{kryterium_2}).
The violation of the RH would manifest as an increase of the amplitude of the combination
$ k^{\frac{3}{4}} c_k$ and for sufficiently large $k$ (depending on $C$) the product
$ k^{\frac{3}{4}} c_k$ will escape outside the strip $\pm C$ (if RH is true we can
take $C=A=7.7751\ldots\times 10^{-5}$). We will discuss below
the case  (\ref{kryterium_2})  of  simple zeros
($\zeta'(s)\neq 0$) requiring $c_k=\mathcal{O}(k^{-3/4})$.

%We will not dwell on
The derivatives $\zeta'(\rho_l)$  in the denominator of (\ref{ck-od-k})  does not pose
any threat to RH. First of all it does not depend on $k$, thus hypothetical extremely
small values of $\zeta'(\rho_l)$ will only change the constant hidden
in big-$\mathcal{O}$ in (\ref{kryterium_2}). Second, this derivative
is taking moderate values for zeros used by us: the smallest $\zeta'(\frac{1}{2}+\gamma_l)$
was  $0.032050162\ldots$  at $\gamma_{1310}==1771.212945\ldots$ and the largest
was $7.7852581838\ldots$ at $\gamma_{1773}=2275.06866478\ldots$.
We remark that Lehman in
\cite{Lehman1960} makes the conjecture:
\bee
\frac{1}{\zeta'(\rho_l)}=\mathcal{O}(\rho_l^\nu),~~~~~~~~{\rm where}~~~~0<\nu<1.
\eee
Some rigorous theorems about the possible large and small values of $\zeta'(\rho_l)$
proved under the assumption of RH can be found in \cite{Ng-2007}. We have checked, that
for first 5549728  nontrivial zeros of $\zeta'(\rho_l)$ the largest derivative was
9.38127677  at $\gamma_{5376610}$ and the smallest was 0.001028760514 at
$\gamma_{4161179}$.

Let us suppose there are some zeros of $\zeta(s)$ off critical line.
Next let us assume that we can split the sum over zeros $\rho_l$ in (\ref{ck-od-k})
in two parts: one over
zeros on critical line and second over zeros off critical line.  This second sum should
violate the overall term $k^{-3/4}$ present in the first sum.  Let $\rho_l^{(o)}$ denote
the zeros lying off critical line (``o'' stands for ``off''):
$\rho_l^{(o)}=\frac{1}{2}\pm\delta_l+i\gamma_l^{(o)},~~\delta_l>0$ (as it is well known the
nontrivial zeros are symmetric with respect to the critical line zero hence the
combination $\pm \delta_l$ plus there are  appropriate conjugate zeros below real axis).
In the factor $|\Gamma(1-\frac{\rho_l}{2})|$ the off-critical line zeros will lead only
to the change of $\gamma_l^\frac{1}{4}  \rightarrow \gamma_l^{\frac{1}{4} \mp \delta_l/2}$  in (\ref{Gamma-exp})
and in view of exponential  decrease  present in (\ref{Gamma-exp}) the possible
violation of RH in (\ref{ck-od-k})
will manifest through the terms $k^{\rho_l/2-1}=k^{-\frac{3}{4}\pm \delta_l/2 +
i\gamma_l^{(o)}/2}$.  The combination $1/k^{\frac{3}{4} + \delta_l/2}$ poses
no problem as it leads to faster than
required in (\ref{kryterium_2}) decrease of some of terms in the sum for $\tilde{c_k}$.
But the combination $1/k^{\frac{3}{4} - \delta_l/2}$ leads to  violation of
(\ref{kryterium_2}) and we want to elucidate how it  can be detected.
Say we want to compare $c_k^g$ with $c_k^e$ with accuracy $\epsilon$, where we
are interested in values of $\epsilon$ of the order $\epsilon=10^{-10\ldots 0}$.
The expression for
$c_k^g$ is a finite sum  and we can in principle calculate its value in PARI with
practically arbitrary exactness (however for really large $k$ it can take years
of CPU time).  Although $c_k^g$ contains information coming from all zeros, to see
influence of the first off critical line zero the value of sufficiently large $k$
has to be examined. The sum for $c_k^e$ is infinite and we
expect that to get accuracy $\epsilon$ we have to sum in (\ref{ck-oscillation})
up to $l=L$ given by (as we skip $|\gamma_l|^{1/4}$ we will skip also $\pi/4$ as our
consideration are not rigorous in general):
\bee
\epsilon\approx e^{-\gamma_L}, {\rm ~~hence}~~  \gamma_L\approx -\log(\epsilon).
\eee
Because values of the imaginary parts $\gamma_l^{(o)}$ of the hypothetical zeros
off critical line should be extremely large, perhaps even as large as $10^{10^{10000}}$,
we suppose that $\gamma_L < \gamma_l^{(o)}$.   The contribution of
$\gamma_l^{(o)}$ is present in $c_k^g$, but will not be present in the explicit
sum  for $c_k^e$  cut at $L$. To detect discrepancy between $c_k^g$ and $c_k^e$
larger than assumed accuracy $\epsilon$ sufficiently large value of $k=K$ is needed.
The point is that $k^{3/4}c_k^g$ will escape outside the strip $\pm C$ for sufficiently
large $k=K$ and the value of such $K$ we can estimate analyzing the explicit formula
for $\tilde{c_k}$.

We can estimate value of  the index $K$ from the requirement  that
the  term $K^{\delta_l/2}$ ($K^{-3/4}$ is present in front of the sum for $\tilde{c_k}$)
  will defeat the
smallness of the term $\Gamma(3/4 - \gamma_l^{(o)}/2)$ and together their
product will overcame the first summand in (\ref{ck-od-k}) corresponding to $\gamma_1$.
In other words in the series (\ref{ck-od-k}) all terms up to $\gamma_l^{(o)}$
monotonically and fast decrease but the terms corresponding to zeros off
critical line can be made arbitrary large for sufficiently large $k$.
The condition for such a $K$ is roughly:
\bee
K^{\delta_l/2} e^{- \gamma_l^{(o)}} >C ~~{\rm thus~~}
K>C'e^{2\gamma_l^{(o)}/\delta_l}.
\eee
Because $\delta_l$ can be arbitrarily close to zero and, as we expect,
$\gamma_l^{(o)}$ is very  large  the value of $K$ will be extremely huge --- larger than
famous Skewes number and will look something like $10^{{10}^{{.}^{{.}^{.}}}}$.
Thus it is practically
impossible to disproof RH by comparing $c_k^g$ and $c_k^e$. \Mas has given in
 \cite{Maslanka3}  discussion  of possible violation of
(\ref{kryterium_2}) and he also came to the pessimistic conclusion that
disproving  RH by comparing $c_k^g$ and $c_k^e$ is ``far beyond any
numerical capabilities'', see pp. 7-8.  We wanted to find agreement between
$c_{100000}^g$ and
$c_{100000}^e$ within one thousand digits  and to our surprise the first attempt
to calculate $c_{100000}^e$ resulted in the difference already on the 87 place.
We started  to struggle with numerical problems  to improve the accuracy and finally
we got 996 digits of  $c_{100000}^g$ and  $c_{100000}^e$  the same.

\section{The computer experiment}

The idea of the experiment is to calculate to high precision the values
of $c_{100000}^g$  and  $c_{100000}^e$ and try to find discrepancy between them.
We calculated one value $c_{100000}^g$ from the generic formula (\ref{ck-generic}),
which contains contribution from all zeros of $\zeta(s)$, even hypothetical zeros
with $\Re(\rho_l)\neq \frac{1}{2}$.
Because $\zeta(2n)$ very quickly tend 1 to get the firm value of $c_k$ it is necessary
to perform calculations with many digits accuracy. Additional problem is fast growing
of binomial symbols. We performed calculation of $c_{100000}^g$ using the free
package PARI/GP \cite{PARI}. This package allows to perform very fast  computations
practically of arbitrary
precision.  We have set precision to 100000 decimals and below in Table I
are the partial sums of (\ref{ck-generic}) recorded after summation of 10000, 20000, ...
100001 terms.  In the middle of  computations the partial sums for $c_{100000}^g$
were of the order $10^{30100}$  to drop finally to $ 1.609757993\ldots\times 10^{-9} $
after summing up all 100001 terms. Separately
we have repeated this calculations with precision set to 150000 places and we
found  the difference only from the 69900 place.  Because the sum (\ref{ck-generic})
is finite we got accurate say  50000 digits from this generic formula for $c_{100000}^g$.
This was an easy part, which took approximately 14 hours for precision 100000
decimals and almost 20 hours for precision 150000  digits on the AMD Opteron 2.6
GHz 64 bits processor.

\vskip 0.4cm
\begin{center}
{\large
{\sf Table I} \\
\bigskip
\begin{tabular}{|c|c|}  \hline
$ n $ & $ \sum_{j=0}^n {(-1)^j \binom{100000}{j}\frac{1}{\zeta(2j+2)}}$ \\ \hline
 $ 10000 $ & $ 5.65168726144550\times 10^{14115} $   \\    \hline
 $ 20000 $ & $ 4.00927204946289\times 10^{21729} $   \\    \hline
 $ 30000 $ & $ 6.08771775660005\times 10^{26526} $   \\    \hline
 $ 40000 $ & $ 5.17938759373151\times 10^{29225} $   \\    \hline
 $ 50000 $ & $ 1.26030418446100\times 10^{30100} $   \\    \hline
 $ 60000 $ & $ 3.45292506248767\times 10^{29225} $   \\    \hline
 $ 70000 $ & $ 2.60902189568574\times 10^{26526} $   \\    \hline
 $ 80000 $ & $ 1.00231801236572\times 10^{21729} $   \\    \hline
 $ 90000 $ & $ 6.27965251271723\times 10^{14114} $   \\    \hline
 $ 100000 $ & $ 1.60975799392038\times 10^{-9} $   \\    \hline
\end{tabular}
}
\end{center}

Now we turn to the calculation
of $c_{100000}^e$ from the explicit formulas (\ref{ck-trend}) and (\ref{ck-oscillation}),
which are infinite sums. It is important not to make the replacement $k+1\rightarrow k$
even for $k=100000$ if we want to get accuracy of the order of 1000 digits.
The series in (\ref{ck-trend}) decreases very fast to zero and
it is very easy to get arbitrary number of digits of  $\bar{c}_k$ using the PARI/GP
procedure \verb"sumalt".

Next we want to calculate $\tilde{c}_{k}$  from  (\ref{ck-oscillation})
with 1000 digits accuracy. From  $e^{-\pi\gamma_l/4}=10^{-1000}$ we get
$\gamma_l\approx 2931.7$ and a glance at the list of zeros of $\zeta(s)$
(e.g. \cite{Odlyzko_zeros}) gives  $l=2402$ because $\gamma_{2402}=2931.0691\ldots$.
First 100 zeros of $\zeta(s)$  accurate to over 1000 decimal places we have taken from
\cite{Odlyzko_zeros}.
Next 2500 zeros of $\zeta(s)$   and derivatives $\zeta'(\rho_l)$  with precision 1000
digits we decided to calculate using the built in Mathematica v.7
procedures \verb"ZetaZero[m]" and numerical differentiation \verb"ND[...]".
%(it was used also by \Mas in )  \verb"ND[Zeta[s]]". \text{1/ND[Zeta[s],s,1/2+I g, Terms-$>$20, WorkingPrecision-$>$1000],;]]}
After a few days we got the values of $\gamma_l$ and $\zeta'(\rho_l)$.
We have checked  using PARI/GP that these zeros were accurate within at least 996
places in the sense that in the worst case $|\zeta(\rho_l)|<10^{-996}, ~l= 1, 2, \ldots, 2600$.
The formulas (\ref{ck-oscillation}) and (\ref{ck-trend}) were implemented in PARI
with precision set to 1000 digits and to our surprise   we obtained  that
the values of $c_{100000}^g$  and $c_{100000}^e$ coincide only up to 87 place:
\bee
10^{-86}>\left|\frac{c_{100000}^g}{c_{100000}^e} -1\right|> 10^{-87}.
\eee

The arguments given earlier in Sect. 3  suggested that no discrepancy should be
found between $c_k^g$ and $c_k^e$  for $k=100000$ hence we  have paid
attention to the numerical inaccuracies in the computation
of the derivative $\zeta'(\rho_l)$ as a possible explanation.
We have played with different options 20, 30, 40 terms and \verb"WorkingPrecision"
in \verb"ND[...]", but finally we got only
the moderate  improvement: the difference between $c_{100000}^g$  and $c_{100000}^e$
shifted  to 105 place:
\bee
10^{-105}< \left|\frac{c_{100000}^g}{c_{100000}^e} -1\right|< 10^{-104}.
\eee

Because in (\ref{ck-oscillation}) gamma functions can be calculated in PARI with
practically arbitrary digits of accuracy the only way to  improve accuracy of
calculation of $c_{100000}^e$ is to find a reliable method of calculating $\zeta'(s)$
with certainty that say all 1000 digits are correct. From (\ref{zeta})
it is easy to obtain the modified expression for zeta:
\bee
\zeta(s)=\frac{1}{1-2^{1-s}}\sum_{n=1}^\infty \frac{(-1)^{n-1}}{n^s}.
\label{zeta-2}
\eee
This series is absolutely convergent for $\Re(s)=\sigma>0$ and can be differentiated
term by term:
\bee
\zeta'(s)=-\frac{\log(2)2^{1-s}}{(1-2^{1-s})^2} \sum_{n=1}^\infty \frac{(-1)^{n-1}}{n^s}+
\frac{1}{1-2^{1-s}}\sum_{n=1}^\infty \frac{(-1)^n\log (n)}{n^s}.
\label{zeta-pochodna}
\eee

\begin{figure}%[pht]
\vspace{-2.5cm}
%\begin{minipage}{15.8cm}
\hspace{-1.0truecm}
\includegraphics[width=15truecm,angle=0, scale=1.1]{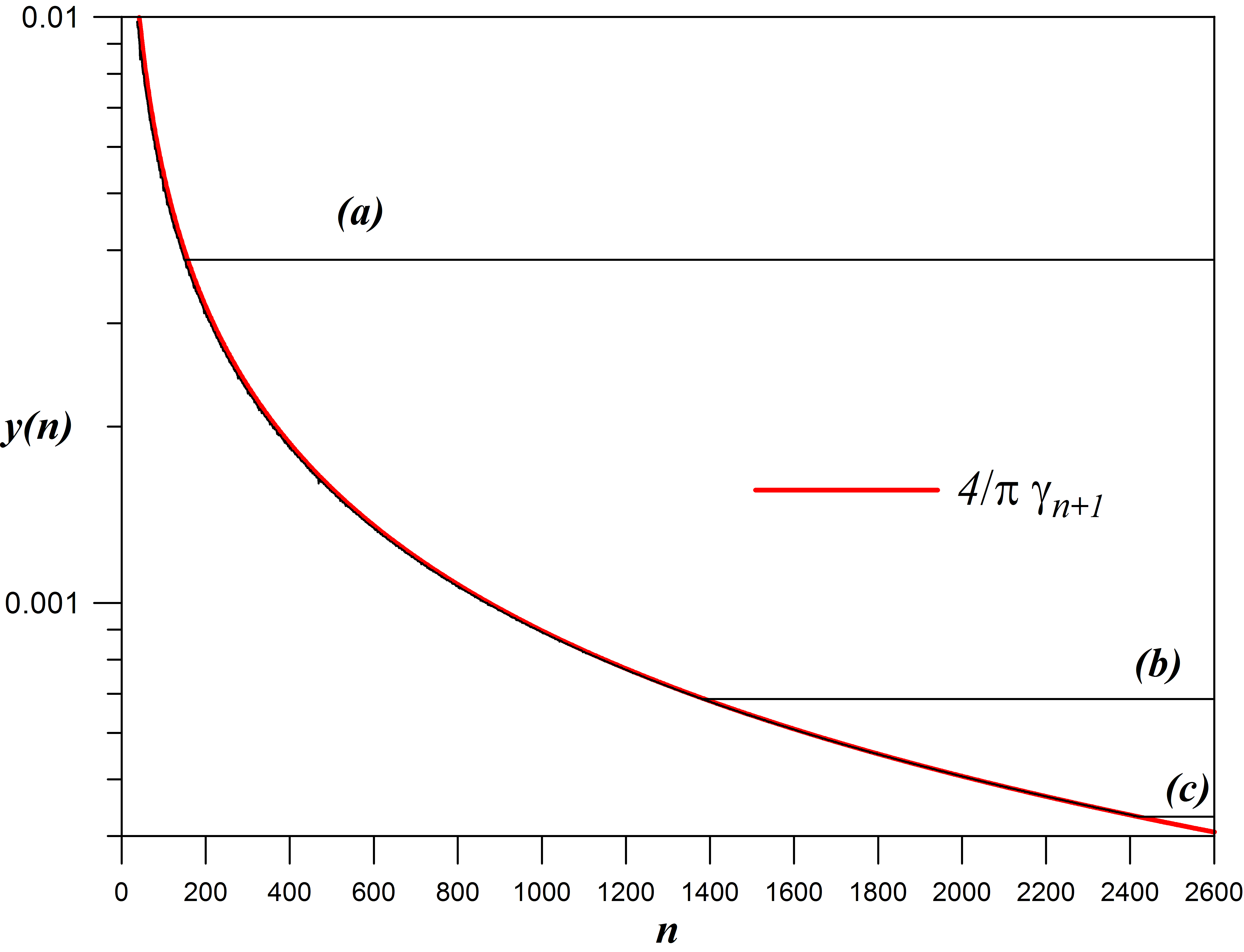} \\
Fig.2  The plot of $y(n)$ for $c_{100000}^e$ obtained by three methods.
The plot $(a)$ is the best result obtained with Mathematica with the procedure
\verb"ND[...]" with  the option \verb"terms=30" and \verb"WorkingPrecision=1000".
The accuracy was 105 digits and the curve $(a)$ departs from the red line at
 the zero $\gamma_{149}=317.73480594237\ldots$ for which
$\exp(-\pi \gamma_{149}/4)=4.193107483\times 10^{-109}$.
The plot $(b)$ was obtained with PARI from the formula (\ref{zeta-pochodna})
using the procedure \verb"sumalt" with precision set to 1000 digits. The $y(n)$ reaches
plateau $0.0006856926750\ldots$ at zero $\gamma_{1412}=1884.00577834967\ldots$, for which
$\exp(-\pi \gamma_{1412}/4)=2.391868726\times 10^{-643}$.
The curve $(c)$ is the same as $(b)$ but derivatives $\zeta'(\frac{1}{2}+i\gamma_l)$
were calculated with \verb"sumalt" and  precision set to 2000 digits. Curve $(c)$ reaches
the plateau $0.00043191361\dots$ at the zero $\gamma_{2430}=2960.033617812\ldots$
for which $\exp(-\pi \gamma_{2430}/4)=2.298783954\times 10^{-1010}$.
\end{figure}

PARI contains the numerical  routine \verb"sumalt" for summing infinite alternating
series  in which extremely efficient algorithm
of Cohen,  Villegas  and Zagier \cite{Zagier-Cohen} is implemented. As this authors
points out on p. 6 their algorithm works even for series like (\ref{zeta-2}) with $s$
complex --- (\ref{zeta-2}) is alternating only when $s\in \mathbb{R}$. We have used
this routine \verb"sumalt" outside scope of its applicability with success
to calculate
$\zeta'(\rho_l)$ from (\ref{zeta-pochodna}) with precision set to 1000 digits and then
calculated $c_{100000}^e$ from
(\ref{ck-oscillation}) and (\ref{ck-trend}). The result was astonishing: the difference
between $c_{100000}^e$ and $c_{100000}^g$ appeared on the 625 place:
\bee
10^{-625}< \left|\frac{c_{100000}^g}{c_{100000}^e} -1\right|< 10^{-624}.
\eee
Because we expect possible violation of RH should manifest at much larger $k$
we  were looking for the way to still improve the accuracy of $\zeta'(\rho_l)$.
We decided to make a frenzy think: we calculated
again $\zeta'(\rho_l)$ using \verb"sumalt"  with zeros having 1000 digits but this
time with precision set to 2000 (however values of $\zeta'(\rho_l)$ were stored only
with 1000 digits).  Thus the aim was to enlarge the number
of terms summed in (\ref{zeta-pochodna}), or in fact the number of iteration performed
inside \verb"sumalt" until the prescribed accuracy is attained.
After 18 hours on AMD Opteron 2.6 GHz we got the results.
And now bingo! The first 996 digits of $c_{100000}^g$ and $c_{100000}^e$ where the same:
\bee
10^{-996}< \left|\frac{c_{100000}^g}{c_{100000}^e} -1\right|< 10^{-995}.
\label{precyzja}
\eee

Because we got the precision (\ref{precyzja}) it is {\it a posteriori} proof that
1000 digits of derivatives $\zeta'(\rho_l)$ were calculated correctly
from (\ref{zeta-pochodna}) by the PARI procedure \verb"sumalt" with precision set
to 2000 digits.

In the Fig.2 we present summary of these computer calculations. Since it is not
possible to plot using standard plotting software as small values as $10^{-600}$ on the
$y$-axis we present in the Fig. 2 the following quantity measuring the distance
from $c_{100000}^g$  to
the partial sums over zeros  $\gamma_l$ in (\ref{ck-oscillation})
and  decreasing with number of  zeros included in the sum:
\bee
y(n)=\left(\log\left(\left|\frac{1}{k^{3/4}}\Re\left(\sum_{l=1}^n
\frac{k^{i\gamma_l/2}\Gamma(\frac{3}{4} -\frac{i}{2}\gamma_l)}
{\zeta'(\frac{1}{2}+i\gamma_l)}\right)-c_k^g\right|\right)\right)^{-1}
\eee
where $k=100000$ and the absolute value is necessary  as the differences between
successive approximants to $c_{100000}^e$ and $c_{100000}^g$ changes sign erratically.
The consecutive terms in the series (\ref{ck-oscillation})
behave like ${\rm e}^{-\pi \gamma_l/4}$ hence we expect that $y(n)$,
by analogy with well known property of alternating series with decreasing terms,
should behave like the first discarded term:
\bee
y(n)\sim \frac{4}{\pi\gamma_{n+1}}.
\label{y-od-n}
\eee
The Fig.2 confirms these considerations: initially $y(n)$ for all curves follows
the prediction (\ref{y-od-n}) and starting with $n$ for which the values of the
$\zeta'(\rho_n)$ are incorrect adding further terms does not improve accuracy.
Because of the exponential decrease of $\Gamma(\frac{3}{4}-\frac{\gamma_l}{2})$
the contribution of further terms is suppressed and horizontal lines in Fig.2 are
determined by the first $\gamma_l$ corresponding to the bad value of the derivative
$\zeta'(\rho_n)$.

\section{Final remarks}

Although we have reached the agreement between $c_{100000}^g$ and $c_{100000}^e$
up to 996 places:
\bdm
c_{100000}^g=1.60975799\ldots~\leftarrow ~980 ~{\rm digits}~ \rightarrow~\ldots 2913696\underline{30140}\times 10^{-9}
\edm
\bdm
c_{100000}^e=1.60975799\ldots~\leftarrow ~980 ~{\rm digits}~ \rightarrow~\ldots 2913696\underline{29833}\times 10^{-9}
\edm
it means nothing about the validity of the RH. The refutation of the RH by computer
methods seems to be as difficult as the analytical proof of its validity. As a possible
precaution let us mention the paper  ``Strange Series and High Precision Fraud''
written by Borwein brothers \cite{Borwein_1992}. In this paper we found a few
striking  examples of the approximate  equalities correct to many, many digits, which
finally are not identities. The most fraudulent is the following
\bee
\sum_{n=1}^\infty \left\lfloor \frac{n e^{\sqrt{163/9}}}{2^n}\right\rfloor \doteq 1280632
\eee
which is valid up to accuracy at least $10^{-500,000,000}$ but is not an identity.

At the end let us remark that the most accurate experiment in physics is the measurement
of the ratio of the electric charge of the electron $e^-$ to the charge of the proton
$e^+$ which is known to be something like $e^-/e^+=-1\pm 10^{-20}$,
see \cite{charges-ratio}. Physicist believe that $e^-=-e^+$ exactly.

\bigskip

{\bf Acknowledgment:} I would like to thank  Dr. Keith Briggs and Prof.
Krzysztof \Mas for  comments and remarks.

\end{document}